\documentclass[10pt]{amsart}
\usepackage{amscd}

\usepackage{amssymb}

\newcommand{\bZ}{{\mathbb Z}}

\newtheorem{thm}{Theorem}[section]
\newtheorem{lemma}[thm]{Lemma}
\numberwithin{equation}{section}

\begin{document}

\title[$\bZ/p^2\times\bZ/p^2$]{Notes on the filtration
of  the $K$-theory for abelian $p$-groups}

\author{Nobuaki Yagita}

\address{ faculty of Education, 
Ibaraki University,
Mito, Ibaraki, Japan}
 
\email{nobuaki.yagita.math@vc.ibaraki.ac.jp, }


\begin{abstract}
Let $p$ be a prime number.  For a given finite group $G$, let $gr_{\gamma}^*(BG)$
 be the associated ring 
of the gamma filtration of the 
$K$-theory for the classifying space $BG$.
In this paper, I correct errors in my paper [Ya]
about $gr_{\gamma}^*(BG)$ when $G$ are abelian 
$p$-groups which are not elementary.
We also extend  Chetard's results for such $2$-groups
to $p$-groups for odd $p$. 
\end{abstract}

\maketitle

\section{Introduction}

Let $p$ be a prime number.  For a given finite group $G$, let $gr_{top}^*(BG)$
(resp. $gr_{\gamma}^*(BG)$) be the associated ring 
of the topological (resp. gamma) filtration of the 
$K$-theory $K^0(BG)$ for the classifying space $BG$.
In this paper, we correct errors in the paper [Ya]
and extend results by B.Chetard.

In Theorem 4.1 in [Ya], I wrote
\begin{thm} ($p=2,r=1$ case by Atiyah [At])
Let $q=p^r$ and $G=\oplus ^n\bZ/q$.  Then
\[ gr_{top}^*(BG)\cong \bZ[y_1,...,y_n]/(qy_i,y_i^qy_j-y_iy_j^q|1\le i,j\le n),\quad |y_i|=2.\]
Hence the filtrations are the same.
\end{thm}

This theorem was error for $r\ge 2$, indeed, arguments for the higher Bokstein
$Q_0'$ in its proof were errors.
However  the statement holds still  for
$p=q$, i.e., for an elementary abelian $p$-group $G$.
(The second statement holds for all abelian $p$
groups [At].)  

Beatrice Chetard pointed out this fact [Ch].  Indeed, she shows the following isomorphism by using  the definition of the gamma filtration
of the representation ring
\[gr_{\gamma}^*(B(\bZ/4\times \bZ/4))\cong
\bZ[y_1,y_2]/(4y_1,4y_2,2y_1^2y_2+2y_1y_2^2,y_1^4y_2^2-
y_1^2y_2^4).\]
She also computes $gr_{\gamma}^*(B(\bZ/4\times
\bZ/2))$, and conjectured
\[gr_{\gamma}^*(B(\bZ/2^r\times \bZ/2))\cong 
\bZ[y_1,y_2]/(2^ry_1,2y_2,y_1y_2^{r+1}+y_1^2y_2^r).\]

In this note, we see that the above Chetard results
can be extended to odd prime cases.
Let us write  $y(1)=y_1^{p}y_2-y_1y_2^{p}$.
 Then  we have
\begin{thm}  For each prime $p$, let $G=\bZ/p^2\times \bZ/p^2$.  
Then  \[  gr_{top}(BG)
\cong \bZ[y_1,y_2]/(p^2y_1, p^2y_2, py(1),y(1)^p).\]
\end{thm}
\begin{thm}
For each prime $p$, let $G=(Z/p^r\times Z/p)$, $r\ge 1$. Then
\[ gr_{top}^*(BG)\cong \bZ[y_1,y_2]/(p^ry_1,py_2,s_r)\]
where $s_r=y_1y_2^{r(p-1)+1}-y_1^py_2^{(r-1)(p-1)+1}$
$mod(y_1^2y_2^{(r-1)(p-1)+2})$.
\end{thm}

Here note that $gr_{\gamma}^*(BG)$ are  known
for many nonabelian $p$-groups $G$ by using $gr_{top}^*(BG)$ and the Atiyah-Hirzebruch spectral sequence, while  the direct computations of $gr_{\gamma}^*(BG)$ by using representations theory 
of $G$ are not so many.

For example, when 
$|G|=p^3$ and nonabelian, we know [Ya]
\[ gr_{\gamma}^*(BG)\cong H^{even}(BG)/(Q_1H^{odd}(BG)).\]
Here $H^{odd}(BG)$  is just $p$-torsion and we can define the Milnor $Q_1$-operation on $H^{odd}(BG)$
(see the proof of Theorem 4.2 in [Ya]). 
In particular, when $G=Q_8$, it is  known $H^{odd}(BG)=0$.
Atiyah [At] used representation arguments to get the multiplicative structure of $H^*(BQ_8)\cong gr_{\gamma}^*(BG)$.

I thank Beatrice Chetard who pointed out my error
in [Ya].

\section{$H^*(B(\bZ/q\times \bZ/q))$ and $H^*(B\bZ/q\times B\bZ/p)$}

Let $X=B\bZ/q$ with $q=p^r$, $r\ge 2$.  Its integral cohomology is
$H^*(X)\cong \bZ[y]/(qy)$ with the degree $|y|=2.$
Considering the long  exact sequence for $q'=q$ or $q'=p$
\[  ...\to H^{*-1}(X;\bZ/q')\stackrel{\delta}{\to} H^*(X)\stackrel{q'}{\to}
H^*(X)\to H^*(X;\bZ/q')\to ...,\]
 we have 
\[ H^*(X;\bZ/q)\cong H^*(X)/q\{1,x\},\quad x=\delta^{-1}y\]
\[ H^*(X;\bZ/p)\cong H^*(X)/p\{1,x'\} \quad x'=\delta^{-1}(p^{r-1}y).\]
Here for a ring $A$, the notation $A\{x,...,z\}$ means
the $A$-free module generated by $x,...,z$.

We consider the Serre spectral sequence
\[ E_2^{*,*'}\cong H^*(X;H^{*'}(X))\Longrightarrow
              H^*(X\times X)\]
\[ where\qquad E_2^{*.*'}\cong \begin{cases}\bZ[y_1]/(qy_1) \quad *'=0\\
\bZ/q[y_1]\{1,x\}\otimes y_2^{*'} \quad  *'>0.
\end{cases}\]
Since $H^*(X)\subset H^*(X\times X)$, element
$y_1^3$ is permanent , and so is $x$.
Hence we have
\[ E_{\infty}^{*,*'}\cong  \bZ/q[y_1,y_2]\{1,y_2x\}.\]
Writing by $\alpha\in H^*(X\times X)$ which 
represents $y_2x\in E_{\infty}^{1,2}$,  we have
\begin{lemma}
For $X=B\bZ/q$, $q=p^r$, we have
\[ H^*(X\times X)\cong \bZ[y_1,y_2]\{1,\alpha\}/(qy_i,q\alpha),\quad |\alpha|=3.\]
\end{lemma}

Next, we consider $X\times B\bZ/p$.  Consider the spectrall sequence
\[ E_2^{*,*'}\cong H^*(X;H^{*'}(B\bZ/p))\Longrightarrow
              H^*(X\times B\bZ/p)\]
\[ where \qquad  E_2^{*.*'}\cong \begin{cases}\bZ/(q)[y_1] \quad *'=0\\
\bZ/p[y_1]\{1,x'\}\otimes y_2^{*'} \quad  *'>0.
\end{cases}\]
\begin{lemma}
For $X=B\bZ/q$, $q=p^r$, we have
\[ H^*(X\times B\bZ/p)\cong \bZ[y_1,y_2]\{1,\alpha'\}/(qy_1,py_2, p\alpha'),\quad |\alpha'|=3.\]
\end{lemma}


\section{$gr_{top}^*(X\times B\bZ/p)$}

In this section, we will prove
\begin{thm}
Let $G=(\bZ/p^r\times \bZ/p)$. Then
\[ gr_{top}^*(BG)\cong \bZ[y_1,y_2]/(p^ry_1,py_2,s_r)\]
where $s_r=y_1y_2^{r(p-1)+1}-y_1^py_2^{(r-1)(p-1)+1}$
$mod(y_1^2y_2^{(r-1)(p-1)+2})$.
\end{thm}

At first, we study relations in $K^*(BG)$.  Recall that
$[p](y)$ is the $p$-th product of the formal group law
of the localized (Morava) $K$-theory ([Ha], [Ra])
so 
that
\[K^*(B\bZ/p)\cong K^*[[y]]/([p](y))\quad |y|=2.\]
We can identify 
$ K^*=\bZ_{(p)}[v_1,v_1^{-1}],$ with $ |v_1|=-2(p-1)$, and 
write
\[[p](y)=py+v_1y^p.\]
We know $K^*(B\bZ/q)\cong K^*[[y]]/([q](y))$ and
\[ [p^r](y)=[p^{r-1}]([p](y))=p^ry+p^{r-1}v_1y^p+... \]

Hence we have  
\[K^*(BG)\cong K^*(X)\otimes_{K^*}K^*(B\bZ/p)
\cong K^*[[y_1,y_2]]/([p^r](y_1),[p](y_2))\]
 \[ where \quad \begin{cases}
[p^r](y_1)=p^ry_1+p^{r-1}v_1y_1^p+...,\\
 [p](y_2)=py_2+v_1y_2^p.
\end{cases} \]
The second equation implies
\[ p^ry_2=-p^{r-1}v_1y_2^p=p^{r-2}v_1^2y_2^{2p-1}=...=
(-1)^rv_1^ry_2^{r(p-1)+1}.\]

Applying this equation to $y_2[p^r](y_1)=0$
in $K^*(BG)$, we get 
\[ y_2[p^r](y_1)=p^ry_1y_2+p^{r-1}v_1y_1^py_2+...\]
\[= (-1)^rv_1^ry_1y_2^{r(p-1)+1}+(-1)^{r-1}v_1^ry_1^py_2^{(r-1)(p-1)+1}+...= (-1)^rv_1^rs_n.\]
Hence $s_n=0$  in $K^*(G)$.
\begin{lemma}  Let $K^*(BG)(y_2)$ be the ideal generated by $y_2$ in  $K^*(BG)$.  Then we have
\[  K^*(BG)(y_2)\cong K^*/p[[y_1,y_2]]\{y_2\}/(s_r).\]
\end{lemma}

We study the Atiyah-Hirzebruch spectral sequence
\[ E_2^{*,*'}\cong H^*(X\times B\bZ/p)\otimes K^{*'}
\Longrightarrow
              K^*(X\times B\bZ/p).\]
Here we recall 
\[H^*(X\times B\bZ/p)\cong
\bZ[y_1,y_2]\{1,\alpha'\}/(qy_1,py_2,p\alpha')\quad  with \
|\alpha'|=3.\]
Since $K^*(X\times B\bZ/p)$ is generated by even dimensional elements, there are $t,t'>1$ and 
$s'\not =0$ in $H^{even}(BG)/p\{y_2\}$ such that
$d_t(\alpha')=v_1^{t'}\otimes s'$.

The map 
\[ v_1^{-t'}\otimes d_t :H^{odd}(BG)\cong \bZ/p[y_1,y_2]\{\alpha'\} \to Z/p[y_1,y_2]\{s'\}\subset \bZ/p[y_1,y_2]\{y_2\}\]
(via $\alpha'\mapsto s'$) is injective. Hence we get 
\[ E_{t+1}^{*.0}\cong  
   \bZ/p[y_1]^+/(qy_1)\oplus \bZ/p[y_1,y_2]\{y_2\}/(s').\]
This term is generated by even dimensional elements,
and isomorphic to
\[ E_{t+1}^{*,0}\cong E_{\infty}^{*,0}
\cong gr_{top}^*(X\times B\bZ/p).\]

Hence $s'$ divides $s_r$. Moreover we can take $s'=s_r$
since  $s'$ $mod(v_1)$  must be in the ideal $(s_r)$
from lemma.

\section{$gr_{top}^*(B\bZ/p^2\times B\bZ/p^2)$}

 Let $X=B\bZ/p^2$.   We study the Atiyah-Hirzebruch spectral sequence
\[ E_2^{*,*'}\cong H^*(X\times X)\otimes K^{*'}
\Longrightarrow
              K^*(X\times X).\]
Here we recall 
$H^*(X\times X)\cong
\bZ[y_1,y_2]\{1,\alpha\}/(p^2y_i,p^2\alpha)$
with $|\alpha|=3.$
We will prove
\[ d_{2p-1}(\alpha)=v_1\otimes py(1),\quad
 d_{2p^2+2p-3}(p\alpha)=v_1^{p+2}\otimes y(1)^p\]
for $y(1)=y_1^py_2-y_1y_2^p$.
Then we see that
\begin{thm}  Let $X=B\bZ/p^2$.  Then we have the isomorphism
\[ gr_{top}^*(X\times X)\cong 
\bZ[y_1,y_2]/(p^2y_i,py(1),y(1)^p).\]
\end{thm}

Recall that the  $p$-product of the formal group law for $K^*$-theory is given by $ [p](y)=py+v_1y^p$.
Then   
\[ [p^2](y)= p(py+v_1y^p)+v_1(py+v_1y^p)^p\]
\[  =p^2y+pv_1y^p+v_1^{1+p}y^{p^2}+...=0\]
in $K^*(X)$.
Consider  in $K^*(X\times X)$, and we have 
\[  y_2[p^2](y_1)-y_1[p^2](y_2) \]
\[= pv_1(y_2y_1^p-y_1y_2^p)+v_1^{p+1}(y_2y_1^{p^2}-y_1y_2^{p^2})+...=0.\]
Let us write $y(i)=y_1^{p^i}y_2-y_1y_2^{p^i}$ so that
\[ (*)\quad py(1)+v_1^py(2)+...=0\quad in \ \  K^*(X\times X).\]

Here we recall the connected $K$-theory $k^*(X)$
such that 
\[k^*=\bZ_{(p)}[v_1],\quad and \quad k^*(X)[v_1^{-1}]=
K^*(X). \]
 Then we can write the above equation $(*)$
\[(**)\quad py(1)+v_1^py(2)=0\ \  mod(p^2v_1y(1))\quad in \ \  k^*(X\times X).\]

Hence  $py(1)$ is zero
in $gr_{top}^*(G)$ but it is nonzero  in $K^*(X\times X)$
(hence $y(2)\not =0$)
since $K^*(X\times X)$ is torsion free.  Therefore we 
see
\[ d_{2p-1}(\alpha)=v_1\otimes py(1),\quad 
and \quad E_{2p}^{*,0}\cong \bZ[y_1,y_2]/(p^2y_i, py(1))
\{1,p\alpha\}. \]

Since $K^*(BG)$ is generated by even dimensional
elements,
we see $\alpha''=p\alpha$ is not a permanent cycle, i.e.
there are  $r>2p, t''>1$, $d\in E_r^{*,0}$ such that $d_r(\alpha'')=v_1^{t''}\otimes d\not =0$.

We study this $d$.  At first $d$ is invariant $mod(p)$
under the action of $GL_2(\bZ/p)$, namely
$d=b$ or $pb$ with $b\in \bZ/p[y_1,y_2]^{GL_2(\bZ/p)}$.  The invariant ring is known as the Dickson algebra
\[ \bZ/p[y_1,y_2]^{GL_2(\bZ/p)}\cong \bZ/p[y(1),y(2)/y(1)],\quad \]
\[where \quad y(2)/y(1)=y_1^{p(p-1)}+y_1^{(p-1)(p-1)}y_2^{p-1}+...+y_2^{p(p-1)}.\]

Consider the restriction to $K^*(X)$
\[ res( y(2)/y(1))=y_1^{p(p-1)}\not =0\in K^*(X)\cong 
K^*[y_1]/([p^2](y_1).\]
Hence we can not take $d=y(2)/y(1)$ neither 
$d=py(2)/y(1)$.

Moreover we still see that $y(2)$ is nonzero.

Therefore when $|d|\le 2(p^2+p)$, we see
$d=y(1)^i$ for $i\le p$.  Here we consider the restriction
to  
\[K^*(X\times X;\bZ/p)\cong K^*/p[y_1,y_2]/(y_1^{p^2},y_2^{p^2}).\]
Hence $d$ is in the $Ideal(y_1^{p^2},y_2^{p^2})$.
 Thus we see that the possibility of the smallest degree
element for $d$ is $y(1)^p$.

We see $y(1)^p=0$  in $gr_{top}^*(X\times X)$
as follows.  Consider in $k^*(X\times X)$
\[y_2^p[p^2](y_1)-y_1^p[p^2](y_2)\]
\[=p^2y(1)+v_1^{p+1}(y_1^{p^2}y_2^p-y_2^{p^2}y_1^p)=
p^2y(1)+v_1^{p+1}y(1)^p\quad mod(p^2v_1y(1)).\]
Here we can see the following equations  with $ mod(p^2v_1y(1))$
 in $k^*(X\times X)$
\[p^2y(1)=p(py(1))=p(v_1^{p+1}y(2))\quad (from \ (**))\]
\[ =
pv_1^{p+1}y(1)\cdot y(2)/y(1)=v_1^{2+2p}y(2)^2/y(1).\]
Continuing this arguments, we see 
\[v_1^{p+1}y(1)^p=0\
mod(v_1^{p+2})\quad in \ \ k^*(X\times X).\]
Thus we can take $d=y(1)^p$.

We see that   the map
\[ \bZ/p[y_1,y_2]\{\alpha'\} \to \bZ/p[y_1,y_2]\{y(1)\} \]
by $\alpha' \mapsto y(1)^p$ is injective.
Hence $E_{2p^2+2p-3}^{*,*'}$ is generated by even dimensional
elements, and is  isomorphic to the
infinite term $E_{\infty}^{*,*'}$.


\begin{thebibliography}{R-W-Y}







\bibitem[At]{At} 
M. Atiyah,
\newblock Characters and the cohomology of finite groups,\
\newblock \emph{Publ. Math. IHES} \textbf{9} (1961), 23-64.

\bibitem[Ch]{Ch}
B. Chetard.
\newblock Graded character rings of finite groups.
\newblock  arXiv:1808.108108 [math.R1], 2018.




\bibitem[Ha]{Ha}
M.Hazewinkel,
\newblock  Formal groups and applications,
\newblock \emph{Pure and Applied Math. 78, Academic Press Inc.} 
\textbf{} (1978), xxii+573pp. 


\bibitem[Ra]{Ra}
D.Ravenel.
\newblock Complex cobordism and stable homotopy groups of spheres. 
\newblock \emph{Pure and Applied Mathematics, 121. Academic Press}
(1986).




\bibitem[Ya3]{YaF}
N.~Yagita.
\newblock  Note on the filtrations of the $K$-theory.
\newblock \emph{Kodai Math. J.}
{\bf 38} (2015), 172-200.


.





\end{thebibliography}
\end{document}